\documentclass[12pt]{amsart}
\usepackage{epsfig}
\newfont{\sheaf}{eusm10 scaled\magstep1}

\newcommand{\ra}{\ensuremath{\rightarrow}}

\def\eea{\end{eqnarray*}}
\def\bea{\begin{eqnarray*}}
\def\Bbb{\bf}

\def\FF{{\Bbb F}}

\def\Ga{{\Gamma}}

\newcommand{\Proof}{{\it Proof. }}
\newcommand{\QED}{{\hfill $Q.E.D.$}}

\newtheorem{teo}{Theorem}[section]
\newtheorem{df}[teo]{Definition}
\newtheorem{lem}[teo]{Lemma}

\newtheorem{oss}[teo]{Remark}
\newtheorem{prop}[teo]{Proposition}

\newcommand{\C}{\ensuremath{\mathbb{C}}}
\newcommand{\R}{\ensuremath{\mathbb{R}}}

\newcommand{\hol}{\ensuremath{\mathcal{O}}}
\newcommand{\HH}{\ensuremath{\mathbb{H}}}
\newcommand{\PP}{\ensuremath{\mathbb{P}}}

\newcommand{\I}{\ensuremath{\mathcal{I}}}
\newcommand{\SSS}{\ensuremath{\mathcal{S}}}

\newcommand{\sC}{{\mathcal C}}

\newcommand{\B}{\ensuremath{\mathbb{B}}}

\begin{document}

\title[A characterization of surfaces uniformized by  the bidisk] {A
characterization of surfaces whose universal cover is the bidisk  }

\author{Fabrizio Catanese, Marco Franciosi}


\date{March 25, 2008}

\begin{abstract}
    We show that the universal cover of a  compact complex  surface  $X$  
  is the bidisk $\HH \times \HH$, or $X$ is biholomorphic to
$\PP^1 \times \PP^1$, if and only if $K_X^2 > 0$ and there exists an
  invertible sheaf $\eta$ such that $\eta^2\cong \hol_X$ and 
    $H^0(X, S^2\Omega^1_X (-K_X) \otimes \eta) \neq 0$.
The two cases are distinguished by   the second
plurigenus,  $P_2(X)\geq 2$ in the former case,  $P_2(X)= 0$ in the latter.
We also discuss related questions.

    \end{abstract}
\maketitle

\section{Introduction}

The beauty of the theory of algebraic curves is deeply related
to the manifold implications of the:

\begin{teo}[Uniformization theorem of Koebe and Poincar\'e]
A connected and simply connected complex curve 
$  \tilde{C} $ is biholomorphic to:    
\[ 
 \tilde{C} \cong 
  \left\{ \begin{array}{ll}  {\PP}^1  & \mbox{ if } \  g=0 \\ 
\C & \mbox{ if } \  g=1 \\
\HH & \mbox{ if } \   g \geq 2 \\  
\end{array}
\right.
\]
\end{teo}  ($\HH$ denotes as usual the
Poincar\'e  upper half-plane 
$\HH= \{ \tau \in \C :  Im (\tau) > 0\}$, but we shall often refer to it 
as the `disk' since it is biholomorphic to
$\{ z\in \C : ||z||< 1 \}$).

Hence  a smooth
(connected) compact  complex   curve  $C$  of genus $g \geq 1$ 
admits a uniformization in the strong sense {\em (iii)} of the following
definition:

\begin{df}\label{unif}
A connected complex space $X$ of complex dimension $n$
admits a {\bf
uniformization} if one of the following conditions hold:
\begin{enumerate}
\item[{\em (i)}]
 there is a connected open set $ \Omega \subset \C^n$
and a surjective holomorphic map $ f \colon \Omega \ra X$
({\bf weak uniformization}); 
\item[{\em (ii)}]
 there is a connected open set $ \Omega \subset \C^n$
and a properly discontinuous group  $ \Ga \subset   Aut (\Omega)$
such that $   \Omega  / \Ga \cong X $
({\bf  Galois uniformization}).
\end{enumerate}

If $X$ is a complex manifold, there are two stronger properties:
\begin{enumerate}
\item[{\em (iii)}]
 there is a connected open set $ \Omega \subset \C^n$
and a surjective holomorphic submersion $ f \colon \Omega \ra X$
({\bf \'etale uniformization});
\item[{\em (iv)}]
 there is a connected open set $ \Omega \subset \C^n$
biholomorphic to the universal cover of $X$
({\bf  strong uniformization}).
\end{enumerate}

\end{df}

Hence the universal cover of a compact complex curve is completely
determined by  its genus; in particular   $\tilde{C} 
\cong \HH$ if and only if  $g \geq 2$, i.e., ``$C$ {\em is of general type}",  
and we get then an isomorphism of $\pi_1 (C)$ with a Fuchsian group  
 $ \Ga \subset \operatorname{Aut} (\HH)  \cong  \operatorname{PSL}(2,
\R)$. 

In higher dimension  the condition that the universal cover
be biholomorphic to a bounded domain $ \Omega$ is quite exceptional;
but still in the Galois \'etale case, where $   \Omega  / \Ga \cong X $ 
and $\Ga$ acts freely with compact quotient, we have, 
if $ \Omega$ is bounded, that the
complex manifold $X$ has ample canonical bundle $K_X$
(see \cite{siegel}), in particular it is a projective manifold
of general type. 

Even more exceptional is the case where the universal cover
is biholomorphic to a bounded symmetric domain $ \Omega$,
or where there is Galois uniformization (ii) of definition \ref{unif})
with source a bounded symmetric domain,
and there is a vast literature on a characterization of these
properties (cf.  \cite{Yau}, \cite{yau1}, \cite{yau}, \cite{Bea}).

The basic result in this direction is S.T.  Yau's uniformization theorem
(explained in \cite{yau1} and \cite{yau}), and for which a very readable exposition is
contained in the first section of \cite{vz},
enphasyzing the role of polystability of the cotangent bundle
for varieties of general type. One would wish nevertheless
for more precise characterizations of the various possible cases.

 For the sake of simplicity, we shall stick here to the case of smooth complex  surfaces,
where the former problem boils down to two very specific
questions.

\hfill\break {\bf Question.} When is the universal cover of  a compact
complex surface
$X$ biholomorphic to the  two dimensional ball 
$ {\B}_2 := \{ z\in \C^2 : || z||
< 1\} $, respectively to  the bidisk $ \HH  \times \HH$ ? 

The first part of this question is fully answered by the
well-known  inequality by Miyaoka and Yau  (cf. \cite{Miy}, \cite{Yau}
\cite{miy}). 
Setting, as usual,  $K_X=$ the canonical divisor, 
$\chi(X): = \chi (\hol_X)$ the holomorphic Euler characteristic  and $P_2(X)=
h^0(X, 2K_X)$ the second plurigenus of $X$, 
 we have the following characterization:

\begin{teo}[Miyaoka-Yau] Let $X$ be a  compact complex  surface. Then $X
\cong {\B}_2 / \Gamma$ (with $\Gamma$ a cocompact  discrete  subgroup of
$\operatorname{Aut} ({\B}_2)$  acting freely on ${\B}_2$)
 if and  only if 
\begin{enumerate}
\item 
$K_X^2 = 9 \chi(S)$; 
\item  the second plurigenus  $P_2(X) >0$.
\end{enumerate}
\end{teo}

The above well known characterization is obtained combining Miyaoka's result
(\cite{miy}) that these two conditions imply the ampleness of  $K_X$, with Yau's
uniformization  result (\cite{Yau}) which uses 
the existence of a K\"ahler Einstein metric; quite  remarkably,  it is
given purely in terms of certain numbers which are either bimeromorphic
or topological invariants.

In the case where $X = \HH  \times \HH/ \Gamma$, with $\Gamma$ a
discrete cocompact subgroup of 
 $\operatorname{Aut} (\HH \times \HH) $ acting freely, one has
 $K_X^2 =8 \chi(X)$. 

But  Moishezon and Teicher in 
 \cite{MT}  showed the existence of  a simply connected 
 surface of general  type  (whence with $P_2(X) >0$)
having $K_X^2 =8 \chi(X)$, 
so that the above conditions are necessary, but not sufficient. 
We observe however that  (and our contribution here is a by-product
of our attempt to answer the latter question) 
 it is still unknown if there exists a 
 surface of general type with   $\chi(X) =1, K_X^2 =8$ 
which is not uniformized by 
$ \HH  \times \HH$. 
 
  The purpose of this
note is to point  out a  precise characterization of  compact complex  surfaces
  whose  universal cover is the bidisk, and of the quadric $\PP^1 \times\PP^1$,
 discussing   whether some
hypotheses can be dispensed with, and to pose an analogous
question in higher dimension.  
Our characterization, which  is of course based on Yau's results, 
relies  on the following crucial

\begin{df}
Let $X$ be a complex manifold of complex dimension $n$.

Then a {\bf special tensor} is a non zero section 
$ 0 \neq \omega \in H^0(X, S^n\Omega^1_X (-K_X) ) $,
while a {\bf semi special tensor} is a non zero section 
$0 \neq \omega \in H^0(X, S^n\Omega^1_X (-K_X) \otimes \eta)$,
 where  $\eta$ is an invertible sheaf 
such that $\eta^2\cong \hol_X$.

We shall say that $X$ admits  a unique semi special tensor if moreover
$ dim (H^0(X, S^n\Omega^1_X (-K_X) \otimes \eta)) = 1.$
\end{df}

In fact, the existence of such tensors is a fundamental property
of manifolds strongly uniformized by the polydisk as we are now going to
see.

Recall that the group of automorphism of $\HH^n$,  
$\operatorname{Aut} (\HH^n) $, is the semidirect product  
 of $(\operatorname{Aut} (\HH))^n$ with the symmetric group
${\mathfrak S}_n$, hence for every subgroup $\Ga$ 
of  $\operatorname{Aut} (\HH^n) $ we have a diagram:

    \[
  \begin{array}{ccccccc}
   1 \ra &  (\operatorname{Aut} (\HH))^n  & \ra & \operatorname{Aut}
(\HH^n) &  \ra {\mathfrak S}_n &  \ra 1\\
   &     \bigcup & & \bigcup &   \bigcup & \\ 
  1 \ra & \Gamma^{0} & \hookrightarrow & \Gamma  & \ra   H & \ra 1 .\\
   \end{array} 
    \]

\begin{prop}\label{nec}
Let $X = \HH^n / \Ga$ be a compact complex manifold whose universal
covering is the  polydisk $ \HH^n$: then $X$ admits a semi special tensor
 and $K_X$ is ample, in particular $K_X^n > 0$.
\end{prop}

\begin{proof}  In $ \HH^n$ take 
coordinates $\{z_1, \dots , z_n\}$  and define 
$$ \tilde{\omega} : = 
\frac{ \operatorname{d}z_1 \otimes \dots \otimes 
\operatorname{d}z_n}{\operatorname{d}z_1\wedge \dots \wedge
 \operatorname{d}z_n}.
 $$  
 Observe that $\tilde{\omega}$ is clearly 
invariant for $(\operatorname{Aut} (\HH))^n$ and for the alternating
subgroup. Let $\eta$ be the 2-torsion invertible sheaf associated to
the signature character of  ${\mathfrak S}_n$ restricted to $H$.
Then clearly $\tilde{\omega}$ descends to a semi special tensor $\omega
\in H^0(X, S^n\Omega^1_X (-K_X) \otimes \eta)$.  

The other assertions are well known (cf. \cite{siegel} and \cite{k=m}).
\end{proof} 

\begin{oss}
We observe that also  $(\PP^1)^n$ admits the following special tensor
$\omega$,   given on $\C^n \subset (\PP^1)^n$ by 
$ \omega : = 
\frac{ \operatorname{d}z_1 \otimes \dots \otimes 
\operatorname{d}z_n}{\operatorname{d}z_1\wedge \dots \wedge
 \operatorname{d}z_n}.
 $
\end{oss}

In dimension two we have then the following

\begin{teo}\label{unibidisk}  Let $X$   be a  compact complex  surface.

Then the following two conditions:
\begin{enumerate}
\item
$X$ admits a semi special tensor;
\item 
$K_X^2 > 0$
    \end{enumerate}

hold  if and only
if either 

 \begin{enumerate}
 \renewcommand\labelenumi{(\roman{enumi})}

\item $ X \cong \PP^1 \times \PP^1$;  or 

\item  $X \cong  \HH  \times \HH / \Gamma$   (where $\Gamma$ is a cocompact 
discrete  subgroup of $\operatorname{Aut} (\HH \times \HH) $
acting freely ). 
\end{enumerate} 
    
\end{teo}

In particular one has the following reformulation of a theorem of S.T. Yau
(theorem 2.5 of \cite{yau}, giving sufficient conditions for  (ii) to hold).

\begin{teo}{\bf (Yau)} \label{yau}
$X$ is strongly uniformized by the bidisk  if and only if
\begin{enumerate}
\item
$X$ admits a semi special tensor; 
\item 
$K_X^2 > 0$;
\item 
 the second plurigenus  $P_2(X)\geq 1$.
    \end{enumerate}

\end{teo}

One can indeed be even more precise:
\begin{teo}\label{sharp} 
$X$ is strongly uniformized by the bidisk  if and
only if
\begin{enumerate}
\item[
{(1*)}] $X$ admits a unique semi special tensor; 
\item[
{(2)}] $K_X^2 > 0$;
\item [
{(3*)}]  the second plurigenus  $P_2(X)\geq 2$. 

    \end{enumerate}
$X$ is biholomorphic to  $ \PP^1 \times \PP^1$ if and only if
(1*), (2) hold and $P_2(X) = 0.$
\end{teo}

It is interesting to see that none  of the above hypotheses can be dispensed 
with.
	
	    	    \begin{oss} The following examples show the existence 
of surfaces which
satisfy two of the  three conditions stated in Thm. \ref{yau},
respectively in Thm. \ref{sharp},
 but  are not uniformized by the bidisk 
\begin{enumerate}
	 \renewcommand\labelenumi{(\roman{enumi})}
	 \item  $\PP^1 \times \PP^1$ satisfies  {\em (1*)} and {\em (2)};  
	 
	  \item A complex torus $X= \C^2 / \Lambda$   satisfies {\em (1) } and {\em (3)}, 
but neither  {\em (1*)} nor  {\em (3*)} (obviously, it does not satisfy
{\em  (2)}); 
	   
	    \item  $X= C_1 \times C_2  $ with $g(C_1) =1$, $g(C_2)=2$ satisfies
{\em  (1*)}  and {\em  (3*)}, but 
	   its universal cover is
	   $ \tilde{X} \cong  \C \times \HH$. 
\end{enumerate} 
   \end{oss}

The most intriguing examples are  provided by 

\begin{prop}\label{elliptic} 
There do exist properly elliptic surfaces $X$ satisfying
\begin{itemize}
\item
{\em (1)} $X$ admits a   special tensor; 
\item 
 {\em (3*)} the second plurigenus  $P_2(X)\geq 2$; 
\item
$ q (X) : = dim ( H^1 (\hol_X )) > 0 $; 
\item
$K_X^2 = 0$; 
\item
$X$ is not birational to a product.
   \end{itemize}
\end{prop}

We would like to pose then the following 

{\bf Question. \ }{\em  Let $X$ be a surface with $ q(X)= 0$
and satisfying {\em (1*)}  and {\em (3*)}: is then $X$ strongly uniformized by the bidisk?}

Concerning the above question, recall the following 
 
  \begin{df}   
$  \Gamma \subset \operatorname{Aut} (\HH^n) $ is said to be reducible
  if there exists $\Gamma^0$ as above 
(i.e., such that 
$\gamma (z_1, ..., z_n) =( \gamma_1(z_1),..., \gamma_n(z_n))$ for every
$\gamma \in \Gamma^0$) and a decomposition 
$\HH^n = \HH^k \times \HH^h$ (with $ h >0$)
 such that the action of $\Gamma^0$  on  $\HH^k $ is discrete.   
\end{df}  

For $n=2$ there are then only two alternatives:  
 
 \begin{oss}
 Let  $\Gamma \subset \operatorname{Aut} (\HH^2) $ be a discrete cocompact
subgroup acting freely 
 and  let $X = \HH^2/ \Gamma $. Then 
 
 \begin{itemize}
 
 \item $\Gamma$ is reducible if and only if $X$ is isogenous to a product
of curves, i.e., there is a finite group $G$ and two curves of genera at least 2
such that 
 $X \cong C_1 \times C_2 / G$. Both cases   $q(X) \neq 0$,
$q(X) = 0$ can occur here. 
 
 \item $\Gamma$ is irreducible and $q(X) = 0$
( this result holds in all dimensions and is a well-known result of Matsushima
 \cite{Ma}).
 \end{itemize} 
 \end{oss} 
 
Let us try to explain the main idea of our main result.
In order to do this, it is important to make the following
\begin{oss}\label{double}
A complex manifold $X$  admits a semi special tensor if and only if
it has an unramified  cover $X'$ of degree at most two which admits a special
tensor.
\end{oss}
\begin{proof}

Assume that we have an 
invertible sheaf $\eta$ such that $\eta^2\cong \hol_X$,
 $\eta \not\cong \hol_X$. Take the corresponding double connected 
  \'{e}tale covering  $ \pi :  X' \ra X$ and observe that
  $$  H^0(X', S^n\Omega^1_{X'} (-K_{X'})) \cong
H^0(X, S^n\Omega^1_{X} (-K_{X})) \oplus H^0(X, S^n\Omega^1_{X}
(-K_{X}) \otimes \eta).$$ 
Whence, there is a  special tensor on $X'$ if and only if there is
a semi special tensor on $X$. 
   \end{proof}

In dimension $ n =2$ things are easier, since the existence of a special tensor 
$\omega$ is equivalent to the existence of a trace free endomorphism
$\epsilon$ of the  tangent bundle of $X$.

Our proof of Theorem \ref{unibidisk}  consists  essentially  in
finding a decomposition of the tangent bundle 
$T_X$ as a direct sum of two  line bundles $L_1 $ and $L_2$,  which are the 
eigenbundles of an invertible
   endomorphism $\epsilon \in \operatorname{End} (T_X)$
(see \S 2 and \S 3 for details), and then applying  the results 
on surfaces with split tangent bundles as given in
\cite{Bea}.

Since the results on manifolds with split tangent bundles 
 hold in dimension $n
\geq 3$, one has a 
characterization of  compact manifolds strongly uniformized by the polydisk
under  a very strong condition on the semi special tensor 
  $\omega \in H^0(X, S^n\Omega^1_X (-K_X) \otimes \eta) $, which essentially
  corresponds to ask for the local splitting of $\omega$ as the product of
$n$ 1-forms which are linearly independent at each point.
 There remains the problem of finding a simpler characterization.

\section{Preliminaries and remarks}\label{preliminari}

\hfill\break {\bf Notation.}  
  $X$  denotes throughout a  compact complex surface.  We use  standard 
notation of algebraic geometry: $\Omega^1_X$  is the cotangent
sheaf,
  $T_X$ is the holomorphic tangent bundle (locally free sheaf), 
$c_1(X)$, $c_2(X)$ are the Chern classes of $X$;
  $K_X$ is the canonical divisor, and $P_n := h^0(X, nK_X)$ is called
the  $n$-th plurigenus, in particular for $n=1$ we have
   the geometric genus of $X$ $p_g(X): = h^0(X, K_X)$, 
while $q : = h^1(X, \hol_X)$ is classically called 
the irregularity of $X$. Finally,  
  $\chi (X) : = \chi (\hol_X) = 1 +p_g -q$ 
is the holomorphic Euler characteristic. 
  
  With a slight abuse of notation,  we do not distinguish between 
invertible sheaves, line
bundles and divisors, while the symbol 
 $\equiv$ denotes  linear equivalence of divisors.

\hfill\break  First of all  let us recall a result of  Beauville which
characterizes  compact complex surfaces whose  universal 
cover is a product  of
two complex curves (cf. \cite[Thm. C]{Bea}). 

\begin{teo}[Beauville] \label{teo:beauville} 
 Let $X $ be a compact  complex surface.
 The tangent bundle
$T_X$ splits as a direct sum of two line bundles
  if and only if either $X$ is a special Hopf surface
or the universal covering space
of
$X$  is a product 
$U\times V$  of two complex curves and the group $\pi_1(X)$ acts
diagonally on 
$U \times V$.
\end{teo}  

Given a direct sum
decomposition of the cotangent bundle $\Omega^1_X \cong L_1\oplus L_2$, 
 Beauville shows   that $(L_1)^2 =( L_2)^2=0$   (cf. \cite[4.1,
4.2]{Bea}) hence  
 \[
 K_X \equiv L_1 +L_2 \hspace{ 1 cm}  c_1(X)^2  = 2\cdot (L_1 \cdot L_2 )
= 2\cdot c_2(X)
 \]  
The last equality corresponds to $K_X^2 = 8 \chi(X)$.

Let us now consider the  bundle
$\operatorname{End} (T_X)$
 of endomorphisms of the tangent bundle.   We can write
$\operatorname{End} (T_X) = \Omega^1_X \otimes T_X $  and  from the
nondegenerate bilinear  map
$$ \Omega^1_X \times \Omega^1_X \longrightarrow \Omega^2_X \cong  K_X$$ we
see that $T_X=  (\Omega^1_X )^{\vee} \cong \Omega^1_X (-K_X)$.  This
exactly means that we  have an isomorphism 
 $\operatorname{End} (T_X)\cong  \Omega^1_X \otimes \Omega^1_X (-K_X).$
 
 Let us see how  this isomorphism
works  in local coordinates $(z_1, z_2)$. I.e.,
 let us see how an element  
 $\frac{ \operatorname{d}z_i \otimes
\operatorname{d}z_j}{\operatorname{d}z_1 \wedge 
 \operatorname{d}z_2}$ 
 in $ \Omega^1_X \otimes \Omega^1_X (-K_X)$  acts on a vector of  the 
form
 $\frac{\partial}{\partial z_h}$.  We have 
 $$\frac{ \operatorname{d}z_i \otimes
\operatorname{d}z_j}{\operatorname{d}z_1 \wedge 
 \operatorname{d}z_2}  \bigl( \frac{\partial}{\partial z_h} \bigr) =
\left\{ \begin{array}{cl} 
 \frac{ \operatorname{d}z_j}{\operatorname{d}z_1 \wedge 
 \operatorname{d}z_2} & \mbox{ if } h=i  \\
0 & \mbox{ if } h\neq i \\
\end{array} \right. 
$$
 In turn,  $\displaystyle{\frac{\operatorname{d}z_j }{\operatorname{d}z_1
\wedge 
 \operatorname{d}z_2} }$ evaluated on  $\operatorname{d}z_k$ gives
$\displaystyle{\frac{\operatorname{d}z_j \wedge 
 \operatorname{d}z_k}{\operatorname{d}z_1 \wedge 
 \operatorname{d}z_2} }$. 

Therefore  a generic element  $ \displaystyle{ \sum_{i,j} a_{ij} 
\frac{ \operatorname{d}z_i \otimes
\operatorname{d}z_j}{\operatorname{d}z_1 \wedge 
 \operatorname{d}z_2}}$
 corresponds to an endomorphism, which,  with respect to the  basis
$\bigl\{ \frac{\partial}{\partial z_1},\frac{\partial}{\partial z_2}
\bigr\}$ is expressed by the matrix 
$$ \begin{pmatrix} -a_{12} & -a_{22} \\ a_{11} & a_{21} \\ 
\end{pmatrix}
$$  In particular for  the symmetric  tensors (i.e., $a_{12}= a_{21}$), 
respectively for the skewsymmetric tensors  (i.e., $a_{12}=-a_{21},
a_{11}=a_{22}=0$) the following isomorphisms hold:
$$ S^2( \Omega^1_X ) (-K_X) \cong \bigg\{  
\begin{pmatrix} -a & -a_{22} \\ a_{11} & a \\  
\end{pmatrix}   \bigg\} \ ; \ \  \hspace{ 0,5 cm} 
 {\bigwedge}^2( \Omega^1_X ) (-K_X) \cong \bigg\{ 
\begin{pmatrix} b & 0 \\ 0 & b \\  
\end{pmatrix}   \bigg\} 
$$ 
We can summarize the above discussion in the following

\begin{lem}\label{split}
If $X$ is a complex surface there is a natural isomorphism
between the sheaf $S^2( \Omega^1_X ) (-K_X)$
and the sheaf of trace zero endomorphisms of the (co)tangent sheaf
$\operatorname{End}^0 (T_X) 
\cong \operatorname{End}^0 (\Omega^1_X)$.

A special tensor $\omega \in H^0 (S^2( \Omega^1_X ) (-K_X))$ 
with nonzero determinant $ det (\omega) \in \C$ yields an
eigenbundle splitting 
$\Omega^1_X\cong L_1 \bigoplus L_2$ of the cotangent
bundle.

If instead $ det (\omega) = 0 \in \C$,  the corresponding endomorphism
$\epsilon$ is nilpotent and yields an exact sequence of sheaves
 \[ 
 0 \ra L \ra \Omega^1_X \ra {\mathcal{ I}}_Z L (-\Delta) \ra 0
  \]
where $L : = ker ( \epsilon)$ is invertible, $\Delta$ is an effective divisor,
and $Z$ is a 0-dimensional subscheme(which is a local complete intersection).

We have in particular 
$K_X \equiv 2L - \Delta$ and
$c_2(X)= {length}(Z) +
 L\cdot (L-\Delta)$.  

\end{lem}

\Proof
We need only to observe that $ det (\omega)$ is a constant,
since  $det (\operatorname{End} (T_X)) =  det (\operatorname{End}
(\Omega^1_X)) \cong \hol_X$.

If $ det (\omega) \neq 0$, there is a constant $c \in \C \setminus \{ 0\}$
such that  $ det (\omega) = c^2$, hence at every point of $X$ the
endomorphism $\epsilon$ corresponding to the special tensor $ \omega$
has two distinct eigenvalues $ \pm c$.

Let $\omega \in H^0(S^2\Omega^1_X (-K_X) )$, $\omega \neq 0$,  be 
 such that $\det (\omega)=0$. 
 Then the corresponding endomorphism $ \epsilon$
 is nilpotent of order 2, and there exists an open nonempty subset  $U
\subseteq X$ such that 
 $\operatorname{Ker}(\epsilon_{|U})  = 
\operatorname{Im}(\epsilon_{|U})$.  
 At a point $p$ where $\operatorname{rank} (\epsilon)=0$, in local
coordinates the endomorphism $\epsilon$ may be expressed by 
 $$ \begin{pmatrix} a & b \\ c & -a \\ 
\end{pmatrix}  \ \ a,b,c \mbox{ regular functions such that } a^2 = - b\cdot c
$$ 
 Let $\delta := \operatorname{G.C.D.} (a,b,c)$. After dividing by
$\delta$, every prime factor 
 of $a$ is either not in $b$, or not in $c$, thus we can write 
 $$ -b= \beta^2 \hspace{ 1 cm}  c=\gamma^2 \hspace{ 1 cm}  a = \beta\cdot
\gamma $$
 Therefore  we obtain 
 $$  \begin{pmatrix} u\\ v \\ 
\end{pmatrix}  \in \operatorname{Ker}{\epsilon} \Longleftrightarrow 
\left\{\begin{array}{l}  a\cdot u + b\cdot v=0 \\ c\cdot u - a\cdot v =0 
\end{array} \right. \Longleftrightarrow \gamma \cdot u - \beta \cdot v =0 
\Longleftrightarrow 
\begin{pmatrix} u\\ v \\ 
\end{pmatrix} = \begin{pmatrix}
\beta \cdot f\\
\gamma\cdot f \\ 
\end{pmatrix}
$$
  and, writing our endomorphism $\epsilon$ as $\epsilon =  \delta \cdot
\alpha$, we have
 $$ \operatorname{Im}(\alpha) = 
 \left\{\begin{array}{l} 
\beta \cdot \gamma \cdot u - \beta^2 \cdot v = \beta \cdot(\gamma \cdot u
- \beta \cdot v)\\
\gamma^2 \cdot u - \gamma \cdot \beta \cdot v = \gamma \cdot (\gamma \cdot
u - \beta \cdot v)
\end{array} \right.
 $$ 
 Let $Z$ be the 0-dimensional scheme defined by $\{ \beta=\gamma=0\}$ and
$\Delta$ be the Cartier divisor defined
 by $\{\delta=0\}$.

 From the above description  we deduce that the kernel of $\epsilon$
is a line bundle $L$ which fits in the following exact sequence:
 \[ 
 0 \ra L \ra \Omega^1_X \ra {\mathcal{ I}}_Z L (-\Delta) \ra 0.
  \]
 Taking the total Chern classes we infer that: $K_X \equiv 2L - \Delta$ as divisors
on $X$ and
$c_2(X)= {length}(Z) +
 L\cdot (L-\Delta)$.  
\qed

\begin{lem}\label{blowup}
Let $X$ be a complex surface and let $X'$ be the blow up of $X$
at a point $p$.
Then a special tensor $\omega'$ on $X'$ induces a special tensor 
$\omega$ on $X$,
and the converse only holds if and only if $\omega$ vanishes at $p$
(in particular, it must hold : $det (\omega) = 0$).
\end{lem}

\Proof
First of all, $\omega'$  induces a special tensor  on $ X \setminus \{p\}$, and
by Hartogs'  theorem the latter extends to a   special tensor 
$\omega$ on $X$.

Conversely, choose local coordinates $ (x,y)$ for $X$ around $p$ 
and take a local chart of the blow up with coordinates $ (x,u)$
where $y=u x$. 
Locally around $p$ we can write
$$
\omega = \frac{a(\operatorname{d}x)^2 + b (\operatorname{d}y)^2+
c (\operatorname{d}x\operatorname{d}y)}{\operatorname{d}x\wedge
\operatorname{d}y} .$$
The pull back $\omega'$ of $\omega$ is given by the following  expression:  
$$
\frac{a(\operatorname{d}x)^2 + b (u \operatorname{d}x +
x \operatorname{d}u)^2 +
c (u\operatorname{d}x+x\operatorname{d}u) \operatorname{d}x}
{x\operatorname{d}x\wedge
\operatorname{d}u} =$$

$$
= \frac{\operatorname{d}x^2
( a+bu^2+cu) + bx^2\operatorname{d}u^2 +
( 2 bu x + c x ) \operatorname{d}x\operatorname{d}u}
{x\operatorname{d}x\wedge
\operatorname{d}u} ,$$
hence $\omega'$ is regular if and only if 
$ \frac{ a+bu^2+cu} {x}$ is a regular function.

This is obvious if $a,b,c$ vanish at $p$, since then their pull back
is divisible by $x$. Assume on the other side that $a,b,c$ are constant:
then we get a rational function which is only regular if
$a = b = c = 0.$

\qed

\begin{lem}\label{rational}
Let $X$ be a compact minimal rational surface admitting a special tensor
$\omega$. Then $ X \cong \PP^1 \times \PP^1$ if $det (\omega ) \neq 0$.
\end{lem}

\Proof
Assume that $X$ is a $\PP^1$ bundle over a curve $B \cong
\PP^1$, i.e., a ruled surface $ \FF _n$ with $ n \geq 0$. Let $\pi \colon X \ra B$ 
the projection. 

By the exact sequence 
$$ 0 \ra \pi^* \Omega^1_B \ra \Omega^1_X \ra \Omega^1_{X|B} \ra 0$$
and since on a general fibre $F$ the subsheaf $\pi^* \Omega^1_B$ is trivial,
while the quotient sheaf $\Omega^1_{X|B}$ is negative, we conclude that any
endomorphism $\epsilon$ carries $\pi^* \Omega^1_B$ to itself.
If it has non zero determinant we can conclude by Theorem \ref{teo:beauville}
that $ X \cong \PP^1 \times \PP^1$. Otherwise, $\epsilon$ is nilpotent
and we have a nonzero element in $ {\rm Hom}(\Omega^1_{X|B}, \pi^*
\Omega^1_B)$.

Since these are invertible sheaves, it suffices to see when
$$ H^0 ( \hol_X ( 2 \pi ^* K_B - K_X )) \neq 0.  $$
But, letting $\Sigma$ be the section with selfintersection
$\Sigma^2 = -n$, our vector space equals  
$ H^0 ( \hol_X ( 2 \Sigma - (n+2) F )) .  $
Intersecting this divisor  with $\Sigma$ we see that (since each time the intersection number
with $\Sigma$ is negative)  $ H^0 ( \hol_X ( 2 \Sigma - (n+2) F )) = H^0 (
\hol_X ( 
\Sigma - (n+2) F ))=  H^0 ( \hol_X (  - (n+2) F )) = 0. $

There remains the case where $X$ is $\PP^2$.

In this case $\epsilon$ must be a nilpotent endomorphism by Theorem
\ref{teo:beauville}, and it cannot vanish at any point by our previous result
on $ \FF _1$. Therefore the rank of  $\epsilon$  equals 1 at each point.
By lemma \ref{split} it follows that there is a divisor $L$ such that $K_X = 2 L$,
a contradiction.

\qed

\section{Proof of Theorems \ref{unibidisk} and \ref{sharp}}

\begin{proof}
If $X$ is strongly uniformized by the bidisk, then $K_X$ is ample,
in particular $K^2_X \geq 1 $ and, since by Castelnuovo's theorem
$ \chi(X) \geq 1$,
by the vanishing theorem of Kodaira and
Mumford it follows that
$ P_2(X) \geq 2$ (see \cite{cm}).

Thus one direction follows from proposition  \ref{nec},
except that we shall show only later that (1*) holds. 

Assume conversely that $(1), (2)$ hold.  
Without loss of generality we may assume by lemma \ref{blowup} that $X$ is
minimal, since $K_X^2$ can only decrease via a blowup and the bigenus is a
birational invariant.

$K^2_X \geq 1 $ implies that either the
surface
$X$ is of general type, or it is a rational surface. In the latter case we conclude by
lemma \ref{rational}.

Observe that the further hypothesis (3)  (obviously implied by (3*))
guarantees that $X$ is of general type.

Thus, from now on, we may assume  that  $X$ is of general type
and, passing to an \'etale double cover if necessary, that
$X$ admits a special tensor.

By the cited Theorem \ref{teo:beauville} of \cite{Bea} it suffices
to find a decomposition of the cotangent bundle 
$\Omega^1_X$ as a direct sum of two  line bundles $L_1 $ and $L_2$. 

The two line bundles $L_1$, $L_2$ will be given as  eigenbundles of a
diagonizable   endomorphism $\epsilon \in \operatorname{End} (\Omega^1_X)$. 

Our previous discussion shows then that it is sufficient to show that any
special tensor cannot yield a nilpotent endomorphism. 

Otherwise, by lemma \ref{split}, we can write 
 $2L \equiv K_X + \Delta$ and then deduce that $L$ is a big divisor since
$\Delta$ is
 effective by construction and $K_X$ is big because  $X$ is of general
type.  This assertion gives the required contradiction since  by 
 the  Bogomolov-Castelnuovo-de Franchis Theorem (cf. \cite{Bog})  for 
   an invertible subsheaf $L$  of  $\Omega^1_X$ it is   
$h^0(X, mL) \leq O(m) $, contradicting the bigness of $L$. 

There remains to show  (1*).  But if   $h^0(X, S^2\Omega^1_X (-K_X) ) \geq 2$
then, given a point $ p \in X$, there is a special tensor which is not invertible in
$p$, hence a special tensor with vanishing determinant, a contradiction.
 
\end{proof} 

\section{Proof of proposition \ref{elliptic}}

In this section we consider surfaces $X$ with bigenus $P_2(X) \geq 2$
(property (3*)), therefore their Kodaira dimension equals 1 or 2,
hence either they are properly (canonically) elliptic, or they are of general type.

Since we took already care of the latter case in the main theorems \ref{unibidisk} 
and \ref{sharp},
we restrict our attention here to the former case, and try to see when does
a properly elliptic surface admit a special tensor (we can reduce to this situation
in view of remark \ref{double}). We can moreover assume that the associated
endomorphism $\epsilon$ is nilpotent by theorem \ref{teo:beauville}.

Again without loss of generality we may assume that $X$ is minimal by
virtue of lemma \ref{blowup}.

	\Proof
Let $X$ be a minimal properly elliptic surface and let $f : X \ra B$
be its (multi)canonical elliptic fibration.
Write any fibre $f^{-1} (p) $ as $F_p = \sum_{i=1}^{h_p} m_i  C_i$
and, setting $ n_p : = G.C.D. (m_i)$, $ F_p = n_p F'_p$,
we say that a fibre is multiple if $n_p > 1$.
By Kodaira's classification (\cite{Kodaira}) of the singular fibres
we know that in this case $ m_i = n_p ,  \forall i.$

Assume that the multiple fibres of the elliptic fibration are $n_1 F_1', \dots
, n_r F_r'$, and consider the divisorial part of the critical locus
$$ \SSS_p : = \sum_{i=1}^{h_p} (m_i-1 )  C_i , \ \  \SSS : = \sum_{p
\in B}\SSS_p $$ 
 so that we have then the exact sequence
$$  0 \ra f^* \Omega^1_B (  \SSS ) \ra 
\Omega^1_X \ra  \I_{\sC}\ \omega_{X|B} \ra 0,$$
where $\sC$ is a 0-dimensional (l.c.i.) subscheme.

For further calculations we separate the divisorial part of the critical locus
as the sum of two disjoint effective divisors,
the multiple fibre contribution and the rest:
$$ \SSS_m  : = \sum_{i=1}^r
(n_i -1 ) F'_i, \  \hat {\SSS} : = \SSS  - \SSS_m .$$

Let us assume that we have a nilpotent endomorphism corresponding to another
exact sequence 
\[ 
 0 \ra L \ra \Omega^1_X \ra {\mathcal{ I}}_Z L (-\Delta) \ra 0,
  \]
in turn determined by a homomorphism 
$$ \epsilon' :  {\mathcal{ I}}_Z L (-\Delta) \ra L,$$
i.e., by a section $$s \in H^0 (\hol_X (\Delta)) = $$
$$=H^0 (\hol_X (2 L - K_X))
= H^0 (S^2(L) ( - K_X)) \subset H^0 (S^2(\Omega^1_X) ( - K_X)) .$$

We observe that, since $ 2 L \equiv K_X + \Delta$, it follows that,
if $F$ is a general fibre, then
$$ L \cdot F = \Delta \cdot F = 0,$$
hence the effective divisor 	$\Delta$ is contained in a finite union of fibres. 

The first candidate to try with is the choice of $ L = L'$, where we set $L' : = f^* \Omega^1_B ( 
\SSS)$.

To this purpose we recall Kodaira's canonical bundle formula:
$$ K_X \equiv  \SSS_m + f^* (\delta) =  \sum_{i=1}^r
(n_i -1 ) F'_i + f^* (\delta), \ deg (\delta) = \chi(X) -2 + 2b,$$
where $b$ is the genus of the base curve $B$.

Then $H^0 (\hol_X (2 L' - K_X))= H^0 (\hol_X (f^* (2 K_B - \delta)  
+ 2 \SSS  - \SSS_m)$, and we search for an effective divisor linearly
equivalent to 
$$f^* (2 K_B - \delta)  
+ 2 \SSS   - \SSS_m =f^* (2 K_B - \delta)  
+ 2 \hat{\SSS}   + \SSS_m .$$

We claim that  $H^0 (\hol_X (2 L' - K_X))= H^0 (\hol_X (f^* (2 K_B - \delta) )$:
it will then  suffice to have examples where  $|2 K_B - \delta | \neq \emptyset.$ 

{\em Proof of the claim}

It suffices to show that 
$f_* \hol_X(2 \hat{\SSS}   + \SSS_m) = \hol_B $.
Since the divisor $2 \hat{\SSS}   + \SSS_m$ is supported on the singular fibres,
and it is effective, we have to show that,
 for each singular fibre $F_p = \sum_{i=1}^{h_p} m_i  C_i$,
neither $2 \hat{\SSS}_p \geq F_p$ nor 
${\SSS_m}_{,p} \geq F_p$.

The latter case is obvious since ${\SSS_m}_{,p} = (n_p - 1) F'_p <  F_p = n_p  F'_p $.

In the former case, $2 \hat{\SSS}_p = \sum_{i=1}^{h_p} 2 (m_i-1)  C_i$,
but it is not possible that $\forall i$ one has $  2 (m_i-1) \geq m_i $,
since there is always an irreducible curve $C_i$ with multiplicity $ m_i = 1$.

\QED for the claim

Assume that the elliptic fibration is not  a product (in this case there is no special
tensor  with vanishing determinant): then the irregularity of $X$ equals
the genus of $B$, whence our divisor on the curve $B$ has
degree equal to $ 2b-2 - (1-b + p_g(X)) = 3 b -3 - p_g. $

Since $\chi (X) \geq 1$,  $p_g : = p_g (X) \geq b$, and there exist
an elliptic surface $X$ with any $ p_g \geq b$ (\cite{qed}).

Since any divisor on $B$ of degree $\geq b$ is effective, it suffices to choose
$  b \leq p_g \leq 2b-3$ and we get a special tensor with trivial determinant,
provided that $b \geq 3$.

Take now a Jacobian elliptic surface in Weierstrass normal form 
$$ Z Y^2 - 4 X^3 - g_2 XZ^2 - g_3 Z^3 = 0,$$
where $ g_2 \in H^0 (\hol_B (  4 M))$, $ g_3 \in H^0 (\hol_B (  6 M))$,
and assume that all the fibres are irreducible.

Then the space of special tensors corresponding to our choice of $L$
corresponds to the vector space $H^0 (\hol_B ( 2 K_B - \delta)) =
H^0 (\hol_B (  K_B - 6 M))$. It suffices to take a hyperelliptic
curve $B$ of genus $ b = 6 h + 1$, and, denoting by $H$ the hyperelliptic
divisor, set $ M : = h H$, so that $K_B - 6 M \equiv 0$
and we have $h^0 (\hol_X (2L-K_X)) =1$. We leave aside for the time
being the question whether the surface $X$ admits a unique special tensor.

\qed    	   


\hfill\break  {\bf Acknowledgements.}

These research   was performed in the realm  of the
   SCHWERPUNKT "Globale Methoden in der komplexen Geometrie", and of the
FORSCHERGRUPPE 790 `Classification of algebraic surfaces and compact complex
manifolds'.

The second author thanks the Universit\"at Bayreuth for its warm
hospitality 
 in the months of november and december 2006 (where the research was
begun) and the DFG  for supporting his  visit.

We would like to thank Eckart Viehweg and Kang Zuo for pointing out some aspects of
Yau's uniformization theorem that we had not properly credited in the first version.

\noindent {\bf Authors' addresses:}

\bigskip

\noindent Prof. Fabrizio Catanese\\ Lehrstuhl Mathematik VIII, 
Universit\"at Bayreuth, NWII\\
   D-95440 Bayreuth, Germany \\  e-mail: Fabrizio.Catanese@uni-bayreuth.de

\bigskip

\noindent Marco Franciosi \\ Dipartimento di Matemativa Applicata "U.
Dini", Universit\`a di Pisa \\
  via Buonarroti 1C, I-56127, Pisa, Italy \\  e-mail:
franciosi@dma.unipi.it


\bigskip


\begin{thebibliography}{99999}

\bibitem[BPV84]{bpv}
W. Barth, C. Peters, A. Van de Ven,  {\em
Compact complex surfaces.}
    {\bf  Ergebnisse der Mathematik und ihrer Grenzgebiete (3).
Springer-Verlag}, Berlin,(1984).

\bibitem[Bea00]{Bea}  A. Beauville, {\em Complex manifolds with split
tangent bundle} Complex analysis and algebraic geometry, de Gruyter,
Berlin, 2000, 61--70.   
\bibitem[Bog77]{Bog} F. A. Bogomolov {\em Families of curves on a surface
of general type} (Russian) 
 {\bf Dokl. Akad. Nauk SSSR,   236} no.  5,   1977 , p.1041--1044. 

\bibitem[Bom73]{cm}
   E. Bombieri,  {\em Canonical models of surfaces of
general type}, {\bf Publ. Math. I.H.E.S., 42} (1973), 173--219.

\bibitem[Cat07]{qed}
F. Catanese,
``{\em Q.E.D. for algebraic varieties}'',
 {\bf Jour. Diff. Geom. 77} no. 1  (2007) 43--75

\bibitem[Kod60]{Kodaira}
K. Kodaira,
{\em  On compact complex analytic surfaces, I},
{\bf  Ann. of Math. 71 } (1960), 111--152.

\bibitem[K-M71]{k=m}
Kodaira, K.  Morrow, J. {\em Complex manifolds}  Holt, Rinehart and Winston, Inc., 
New-York-Montreal, 
Que.-London, 1971


\bibitem[Ma62]{Ma} Y . Matsushima,  {\em On the first Betti number of
compact quotient spaces of higher-dimensional symmetric spaces} 
 {\bf Ann. of
Math. (2) 75},  1962,  312--330  
\bibitem[MT87]{MT}  B.Moishezon, M. Teicher, {\em Simply-connected
algebraic surfaces of positive index}
 {\bf  Invent. Math. 89}  (1987), no. 3, 601--643. 

\bibitem[Miy77]{Miy} Y. Miyaoka,   {\em  On the Chern numbers of surfaces
of general type} {\bf   Invent. Math. 42}  (1977), 225--237. 

\bibitem[Miy82]{miy} Y. Miyaoka,   {\em  Algebraic surfaces with positive indices}  Classification of algebraic and analytic manifolds (Katata, 1982),   {\bf Progr. Math. 39}  BirkhŠuser Boston, Boston, MA, (1983) 281--301

\bibitem[Mu79]{Mu} D. Mumford, {\em An algebraic surface with $K$ ample,
$K^2=9$, $p_g = q = 0$}
 {\bf  Amer. J. Math. 101}, no. 1,  (1979), 233Ð244. 

\bibitem[Sieg73]{siegel}
C.L. Siegel,
{\em Topics in complex function theory. Vol. III. Abelian functions and
modular functions of several variables. }
Translated from the German by E. Gottschling
and M. Tretkoff. With a preface by Wilhelm Magnus. Reprint of the
1973 original.
{\bf Wiley Classics Library. A Wiley-Interscience Publication. John
Wiley and Sons,
Inc.}, New York (1989), x+244 pp.

\bibitem[V-Z05]{vz}
E.  Viehweg, K. Zuo,
    {\em Arakelov inequalities and the uniformization of certain rigid Shimura varieties},
arXiv:math/0503339. 


\bibitem[Yau77]{Yau} S. T. Yau,  {\em Calabi' s conjecture and some new
results in algebraic geometry}
 {\bf   Proc. Nat. Acad. Sci. USA 74}  (1977),  1798--1799

\bibitem[Yau88]{yau1} S. T. Yau,  {\em Uniformization of geometric structures.}  The mathematical
heritage of Hermann Weyl (Durham, NC, 1987), 
{\bf Proc. Sympos. Pure Math., 48, Amer. Math.
Soc., Providence, RI} (1988),  265--274. 

\bibitem[Yau93]{yau} S. T. Yau,  {\em  
A splitting theorem and an algebraic geometric characterization of locally hermitian symmetric spaces}
{\bf  Comm. in Analysis and Geometry 1}   (1993), 473--486
 
 
   


\end{thebibliography}
\end{document}